\documentclass{amsart}
\usepackage{mathrsfs,eucal,amsmath,amsthm,amsfonts,amssymb,amscd,amsbsy,latexsym}
\usepackage{verbatim}
\usepackage[all,2cell]{xy} \UseAllTwocells \SilentMatrices
\begin{document}


\newtheorem{thm}{\textsc{Theorem}}[section]

\newtheorem{dfn}[thm]{\textsc{Definition}}

\newtheorem{lem}[thm]{\textsc{Lemma}}

\newtheorem{conj}[thm]{\textsc{Conjecture}}

\newtheorem{cor}[thm]{\textsc{Corollary}}

\newtheorem{pro}[thm]{\textsc{Proposition}}
\newtheorem{rem}{\textsc{Remark}}


\newcommand{\matriz}[4]{\left(\begin{array}{cc}#1 & #2 \\ #3 & #4\end{array}\right)}

\newcommand{\aro}{\longrightarrow}

\newcommand{\quadrado}[4]{\xymatrix{#1 \ar[r]\ar[d]&  #2 \ar[d] \\   #4\ar[r] &  #3}}

\newcommand{\quadradod}[4]{\xymatrix{#1 &  #2 \ar[l] \\   #4\ar[u] &  #3\ar[l]\ar[u]}}

\newcommand{\cart}{\ar@{}[dr]|{\square}}

\newcommand{\cartt}{\ar@{}[drr]|{\square}}

\newcommand{\rela}[4]{\xymatrix{#1  \ar@<+.7ex>[r]^{#2}\ar@<-.7ex>[r]_{#3} & #4}}

\newcommand{\seta}[3]{\xymatrix{#1\ar[r]^{#2} & #3}}

\newcommand{\setaa}[3]{\xymatrix{#1\ar[rr]^{#2}&& #3}}


\newcommand{\NN}{\mathbb{N}}

\newcommand{\PP}{\mathbb{P}}

\newcommand{\ZZ}{\mathbb{Z}}

\newcommand{\CC}{\mathbb{C}}

\newcommand{\QQ}{\mathbb{Q}}

\newcommand{\RR}{\mathbb{R}}

\newcommand{\FF}{\mathbb{F}}

\newcommand{\DD}{\mathbb{D}}

\newcommand{\qq}{\overline{\QQ}}


\newcommand{\cc}[1]{\mathscr{#1}}

\title{On the rank of the fibres of rational elliptic surfaces}
\author{Cec\'ilia Salgado}
\address{Niels Bohrweg 1, 2333 - CA Leiden
The Netherlands}



\begin{abstract}
We consider an elliptic surface $\pi: \cc{E}\rightarrow \PP^1$ defined over a number field $k$ and study the
problem of comparing the rank of the special fibres over $k$ with that of the
generic fibre over $k(\PP^1)$. We prove, for a large class of rational elliptic
surfaces, the existence of infinitely many fibres with
rank at least equal to the generic rank plus two. 
\end{abstract}

\maketitle

\section{Introduction}

 Let $k$ be a number field, $B$ a projective curve and $\pi:\cc{E}\rightarrow B$ an elliptic surface over $k$, i.e., a projective surface endowed with a morphism $\pi$ such that almost all fibres are genus-one curves and such that there is a
section $\sigma$ of $\pi$ defined over $k$ that will be fixed as the zero section. The generic fibre $\cc{E}_{\eta}$ is an
elliptic curve over the function field $k(B)$. Since $k(B)$ is finitely generated over $\mathbb{Q}$, the Mordell-Weil Theorem
is still valid in this context (Lang-N{\'e}ron), and hence the set of $k(B)$-rational points of $\cc{E}_{\eta}$ is a finitely generated
abelian group. Since all but finitely many fibres of $\pi$ are elliptic curves over the number field $k$, it is natural to be interested in comparing
the rank $r$ of $\cc{E}_\eta(k(B))$ and the rank $r_t$ of the Mordell-Weil group of a fibre $E_t(k)$ for $t\in B(k)$.

 A theorem on specializations by N\'eron \cite{N} or its refinement by Silverman \cite[Theorem III.11.4]{Sil2}, \cite{Sil} in the case where the base is a curve tells us that if $\cc{E}$ is non-split, then for all but finitely many fibres we have \[r_t\geq r.\]

 Billard in \cite[Theorem C]{B} showed that if we assume $\cc{E}$ to be
 $\QQ$-rational (birational to $\PP^2$ over $\QQ$) and non-isotrivial,
 then \[ \# \{t\in B(k); r_t\geq r+1\}= \infty .\]

 Three natural questions arise:

\begin{itemize}
\item[1)] Can we replace $\QQ$ by an arbitrary number field $k$ and the hypothesis $\cc{E}$ is $k$-rational by $\cc{E}$ is $k$-unirational?
\item[2)] Can we improve the bound, i.e., have $r_t\geq r+2$ for infinitely many $t \in B(k)$?
\item[3)] Can we obtain a similar result for elliptic surfaces that are geometrically rational (i.e., fixed an algebraic closure $\bar{k}$ of $k$, $\cc{E}$ is $\bar{k}$-rational) but not $k$-rational or for non-(geometrically) rational elliptic surfaces, such as K3 surfaces?
\end{itemize}

We remind the reader that a geometrically irreducible algebraic variety $V$ is said to be $k$-unirational if there is a rational map of finite degree $\mathbb{P}^n\dashrightarrow V$ defined over $k$.
 
 We will give a positive answer to question 1) and a partial answer to questions 2) and 3) in this
 article. Question 3) for K3 surfaces is addressed in the author's Ph.D. thesis \cite{salgado} and will
 be explored in another article.

 Let $X$ be a smooth projective rational surface. We denote by $\omega_X^2$ the self-intersection number of the canonical sheaf on $X$. This will be called the degree of $X$ and will be denoted by $d_X$ or by $d$ when the dependency on $X$ is clear. The corollary to the following theorem answers question 1).

\begin{thm}\label{th2}
Let $\pi: \cc{E}\rightarrow B \simeq \mathbb{P}^1$ be a $k$-unirational elliptic surface defined over
a number field $k$. There is a curve $C \rightarrow B$ such that $C \simeq_{k} \mathbb{P}^1$ and \[
\mathrm{rank } (\cc{E}_C(k(C)))\geq \mathrm{rank } (\cc{E}(k(B)))+1,\]
where $\cc{E}_C= \cc{E}\times_B C$.  
\end{thm}

Since the curve $C$ has infinitely many $k$-rational points, an application of N\'eron's or Silverman's Specialization Theorem yields the following corollary.

\begin{cor}\label{cor2}
Let $\pi: \cc{E} \rightarrow B$ be a $k$-unirational elliptic surface then
\[ \# \{ t \in B(k); r_t \geq r+1\}= \infty .\]
\end{cor}

 \begin{rem}\emph{ Since rational surfaces of degree $d_X\geq 3$ such that $X(k)\neq \emptyset$, and del Pezzo surfaces of degree 2 having a $k$-rational point outside a certain divisor, are always $k$-unirational we conclude, from the remark above, that the class of rational elliptic surfaces to which Theorem \ref{th2} applies is quite large. For example, it contains all rational elliptic surfaces defined over $k$ with three distinct types of reducible fibre and/or a fibre with a double component not of type $I_0^*$. But surfaces with generic rank, over $k$, zero and no reducible fibres always have as $k$-minimal models del Pezzo surfaces of degree one and are therefore excluded from the hypothesis of Theorem \ref{th2}.}
\end{rem}

\begin{rem}\emph{ Showing a result such as the above corollary for a rational elliptic surface having as a $k$-minimal model a del Pezzo surface of degree
    one is equivalent to showing that the $k$-rational points on $X$ are Zariski dense, a well known open problem. In
  \cite{Ulas} one can find partial results towards that direction.}
\end{rem}

  The following theorems, or more precisely their corollary, answer question 2) after strengthening the hypothesis of Theorem \ref{th2}.  In order to state them we introduce the following terminology: if $f:X \rightarrow Y$ is a birational morphism of surfaces, passing to an algebraic closure, it is composed of monoidal transformations or blow ups of points; we will refer to the set of $\bar{k}$-points where $f^{-1}$ is not defined as the \textit{blow up locus of} $f$. We will denote it by $\mathrm{Bl}(f)$. Note that if $X$ is a rational elliptic surface and $Y$ a rational model obtained after contracting $(-1)$-curves, then $\mathrm{Bl}(f)$ contains $d_Y$, not necessarily distinct points, thus $\# \mathrm{Bl}(f)\leq d_Y$ where $\#$ denotes the number of distinct points. If $f$ is defined over $k$ the set $\mathrm{Bl}(f)$ is composed of $\mathrm{Gal}(\bar{k}|k)$-orbits.
  
\begin{thm}\label{th1} Let $\pi: \cc{E}\rightarrow B\simeq \PP^1$ be a rational elliptic surface defined over a number field $k$ such that there is a $k$-birational morphism $f: \cc{E} \rightarrow \mathbb{P}^2$ in which the zero section of $\cc{E}$ is contracted to a point $p_1\in \mathbb{P}^2(k)$. Suppose that the blow up locus of $f$ contains at least one orbit distinct from the one given by $p_1$ whose points are, together with $p_1$, in general position. Suppose also that $\cc{E}$ has at most one non-reduced fibre.

Then there exists a finite covering $C\ \rightarrow B$ such that $C(k)$ is infinite and the surface $\cc{E}_C= \cc{E}\times_B C$ satisfies: 
\[ \mathrm{rank}(\cc{E}_C(k(C)))\geq \mathrm{rank}(\cc{E}(k(B)))+2.\] 
\end{thm}

\begin{rem} \emph{It is simple to construct examples of rational elliptic surfaces sa\-tisfying the hypothesis of Theorem \ref{th1}. Let $f$ and $g$ be two arbitrary cubics in $\mathbb{P}^2$, whose equations have coefficients in $k$. Suppose that they pass through two $k$-rational points $p_1$ and $p_2$, then the elliptic surface given by the blow up of the intersection locus of $f$ and $g$ is certainly in this class. }
\end{rem}

\begin{rem}\emph{The assumption that $\cc{E}$ has at most one non-reduced fibre excludes only one configuration of singular fibres, namely, $(I_0^*, I_0^*)$. This case is left out because the surface might become trivial after a quadratic base change (see Lemma \ref{stillrational})}.
\end{rem}

 The result above also holds for some rational elliptic surfaces such that $k$-minimal models, after contracting the zero section, are not isomorphic to $\mathbb{P}^2$ but to other rational surfaces defined over $k$. These are the subject of the next theorem.

 \begin{thm}\label{th3}Let $\pi: \cc{E}\rightarrow B \simeq \PP^1$ be a rational elliptic surface defined over a number field $k$. Suppose $\cc{E}$ does not have reducible fibres. Let $X$ be a $k$-minimal model of $\cc{E}$ such that there exists a birational morphism $f:\cc{E}\rightarrow X$ in which the zero section of $\cc{E}$ is contracted. Suppose that $X$ has degree $d$ and satisfies one of the following:
\begin{itemize}
\item[i)] $d=4,5$ or 8.
\item[ii)] $d=6$,  $\#\mathrm{Bl}(f)=6$ and the largest $\mathrm{Gal}(\bar{k}|k)$-orbit in it has at most four points.
\end{itemize}
Then there exists a finite covering $C\rightarrow B$ such that $C(k)$ is infinite and the surface $\cc{E}_C= \cc{E}\times_B C$ satisfies: 
\[ \mathrm{rank}(\cc{E}_C(k(C)))\geq \mathrm{rank}(\cc{E}(k(B)))+2.\] 
\end{thm}

 
 
Elliptic surfaces satisfying the hypothesis of the previous theorem are always the blow up of a del Pezzo surface $Y$. Indeed, let $Y$ be the surface obtained after contracting the zero section of $\cc{E}$. Since $Y$ is a rational surface, all we have to check is that $-K_Y$ is ample. In fact, the anti-canonical divisor of $Y$ satisfies $(-K_Y)^2>0$ and $-K_Y.D>0$. The former because $K_{\cc{E}}^2=0$ and $Y$ is obtained by contracting a curve in $\cc{E}$. The latter follows from the fact that $\cc{E}$ has no reducible fibres and thus $-K_{\cc{E}}.C \geq 0$ for all $C \in \mathrm{Div}(\cc{E})$ with equality if and only if $C \equiv -K_{\cc{E}}$. By the Nakai-Moishezon Theorem $-K_Y$ is ample.

  Once again an application of N\'{e}ron-Silverman's Specialization Theorem yields the following corollary to Theorems \ref{th1} and \ref{th3}.

\begin{cor}\label{cor1}
Let $\pi: \cc{E} \rightarrow B$ be an elliptic surface as in Theorem \ref{th1} or \ref{th3}. For $t\in B(k)$, let $r_t$ be the rank of the fibre above the point $t$ and $r$ the generic rank. Then \[ \# \{ t \in B(k); r_t \geq r+2\}= \infty .\]
\end{cor}

 \begin{rem}\emph{ Since geometrically, i.e., over $\bar{k}$, a rational elliptic surface is isomorphic to the blow up of nine, not necessarily distinct, points in $\mathbb{P}^2$ (see Proposition \ref{mir}), a $k$-minimal model $X$ of a rational elliptic surface satisfies $1\leq \omega_X^2 \leq 9$. As we suppose that the elliptic surface always has a section defined over the base field $k$ which is contractible, $X$ also verifies $X(k) \neq \emptyset$. }
\end{rem}

\begin{rem} \emph{Theorems \ref{th1} and \ref{th3} are valid for a larger class of rational
  elliptic surfaces. The choice of the cases stated was made for the
  sake of simplicity.  The reader is invited to consult the appendix for examples of cases to which the conclusion of Theorems \ref{th1} and \ref{th3} still applies. } 
\end{rem}

  This text is divided as follows: Sections 2 and 3 contain geometric and arithmetic preliminaries, respectively. Section 4 is dedicated to the proof of Theorem \ref{th2}. It contains a key proposition that reduces the proof of this theorem to the construction of a linear pencil of genus zero curves defined over the base field $k$. The proofs of Theorems \ref{th1} and \ref{th3} are given in Section 5, where we give a case-by-case construction of two pencils of curves of genus zero satisfying certain geometric conditions. The last section sheds some light from analytic number theory into the problem. There, we combine the results obtained in this article with analytic conjectures to get better, but conditional, bounds for the ranks. 
 
\section{Geometric preliminaries}
\subsection{Base change}\label{bc}

 Let $\pi: \cc{E} \rightarrow B$ be an elliptic surface endowed with a section and $\iota: C' \subset \cc{E}$ an irreducible curve. Let $\nu: C \rightarrow C'$ be its normalization. If $C'$ is not contained in a fibre, the composition $\varphi = \pi \circ \iota \circ \nu$ is a finite covering $\varphi: C \rightarrow B$.
 \[
\xymatrix{ & & \cc{E} \ar[d]^{\pi}  \\ C \ar[r]^{\nu} \ar@/_2pc/[rr]^{\varphi} & C'
  \ar[r] \ar[ur]^{\iota} & B } 
\]
We obtain a new elliptic surface by taking the following fibred product:
\[\pi_C: \cc{E}_C= \cc{E}\times_B C \rightarrow C .\]
Each section $\sigma$ of $\cc{E}$ naturally induces a section in $\cc{E}_C$:
 $$(\sigma,id): C= B\times_B C \rightarrow \cc{E} \times_B C.$$ 
 
 We call these sections $\textit{old sections}$. The surface $\cc{E}_C$ has also a \textit{new section} given by 
 $$\sigma^{\text{new}}_C= (\iota \circ \nu, id): C \rightarrow \cc{E}
\times_B C.$$ 
 
 \begin{rem}\emph{Given a smooth elliptic surface $\cc{E}$. The elliptic surface obtained after a base change of $\cc{E}$ is not necessarily smooth. When this is the case, we will replace the base changed surface by its relatively minimal model without further notice. }
\end{rem}

 \begin{rem}
 \emph{ If $C$ is not contained in a fibre nor in a section, then  the new section is different from the \textit{old} ones, but it is not necessarily linearly independent in the Mordell-Weil group. }
  \end{rem}
  
 \begin{rem} \emph{Since we want $B(k)$ to be infinite we are naturally led to consider the base curves $B$ such that $B \simeq \PP^1$ or with geometric genus $g(B)=1$ and $\mathrm{rank}(B(k))\geq 1$.}
 \end{rem}

 If $\cc{E}$ is a rational elliptic surface then, in general, after a quadratic base change we obtain an elliptic K3 surface $\cc{E}'$; nevertheless, if $\cc{E}$ has a non-reduced fibre i.e., a fibre of type $*$, and the base change is ramified above the place corresponding to this fibre and above the place corresponding to a reduced fibre, then the base-changed surface $\cc{E}'$ is still rational. 
\begin{lem}\label{stillrational} Let $\cc{E} \rightarrow B$ be a rational elliptic surface. Let $\varphi: C\rightarrow B$ be a degree two morphism where $C$ is rational. Then one of the following occurs:
\begin{itemize}
\item[i)] $\cc{E}$ has a non-reduced fibre, i.e., a fibre of type $*$, and the morphism $\varphi$ is ramified above the place corresponding to it and above the place corresponding to a reduced fibre. In this case, $\cc{E}_C= \cc{E} \times_B C$ is a rational elliptic surface.
\item[ii)] $\cc{E}$ has two non-reduced fibres, which are necessarily of type $I_0^*$, and the morphism $\varphi$ is ramified above both non-reduced fibres. In this case $\cc{E}_C \simeq E \times D$ where $E$ is an elliptic curve and $D$ is a curve of genus zero. 
\item[iii)]  $\cc{E}_C$ is a K3 surface.
\end{itemize}
\end{lem}
\textbf{}

\noindent \textbf{Proof:} Let $v$ be a place of $B$ and $w$ be a place of $C$ above it. If $w$ is not ramified, then the type of the fibre $(\cc{E}_C)_w$ is the same as that of $(\cc{E})_v$. In this case, since the degree of $\varphi$ is two, there are two places $w, w'$ above $v$. If $w$ ramifies above $v$ and $\cc{E}_v$ is a singular fibre, then fibre type changes, namely:
\begin{itemize}
\item[a)] A fibre $(\cc{E})_v$ of type $*$ induces a fibre $(\cc{E}_C)_w$ whose contribution is $d_w=2d_v-12$ where $d_v$ denotes the contribution of $(\cc{E})_v$, i.e., the local Euler number of $\cc{E}_v$.
\item[b)] A reduced fibre, i.e., a fibre $(\cc{E})_v$ of type $I_n, II, III$ or $IV$ transforms to a fibre $(\cc{E}_C)_w$ whose contribution to the Euler number of the surface is $d_w=2d_v$.
\end{itemize}
 
 Since $\cc{E}$ is rational we have \[\sum_{v \in B} d_v=12.\] \indent
 Thus if $\cc{E}$ has a fibre, $(\cc{E})_v$,  of type $*$ and $\varphi$ is ramified above the place corresponding to it and above the place corresponding to a reduced fibre, then the Euler number of $\cc{E}_C$ is $\sum_{w\in C} d_w= 2d_v-12+2(12-d_v)=12$ implying that $\cc{E}_C$ is rational.

 If $\cc{E}$ has two non-reduced fibres $\cc{E}_{v_1}$ and $\cc{E}_{v_2}$, then since each contributes with at least 6 to the Euler number, which is 12, we have that the contribution of each must be exactly 6. This implies that each non-reduced fibre is of type $I_0^*$ and moreover, that these are the only singular fibres of $\cc{E}$. If $\varphi$ ramifies above both $\cc{E}_{v_1}$ and $\cc{E}_{v_2}$ then by a), they both become non-singular fibres. Since these were the unique singular fibres, the base-changed surface $\cc{E}_C$ has only smooth fibres. This entails $\cc{E}_C\simeq E\times D$ where $E$ is an elliptic curve and $d$ is a curve of genus zero.

Otherwise, i.e., if $\varphi$ is ramified only above reduced fibres, then the Euler number of $\cc{E}_C$ is $2d_v+2(12-d_v)=24$, and hence $\cc{E}_C$ is a K3 surface.

 \subsection{Construction of rational elliptic surfaces}
 
 Let $F$ and $G$ be two distinct cubic curves in $\PP^2$. We will also denote by $F$ and $G$ the two homogeneous cubic polynomials associated to these curves. Suppose $F$ is smooth.
 
  The pencil of cubics generated by $F$ and $G$,
   \[\Gamma:= \{tF+uG: (t:u) \in \mathbb{P}^1\},\] 
   has nine base points (counted with multiplicities) namely the intersection points of the curves $F$ and $G$. The blow up of these points in $\PP^2$ defines a rational elliptic surface $\cc{E}_{\Gamma}$.
   
   Conversely, over an algebraically closed field, we have the following proposition.
   
   \begin{pro}\label{mir}\cite{Mi} Over $\bar{k}$ every rational elliptic surface with a section is isomorphic to a surface $\cc{E}_{\Gamma}$ for a pencil of cubics $\Gamma$ as above.
   \end{pro}
   
   This proposition motivates the following definition.
   
   \begin{dfn}\emph{ Let $\cc{E}$ be a rational elliptic surface. We say that a cubic pencil $\Gamma$ in $\PP^2$ induces $\cc{E}$ if $\cc{E}$ is $\bar{k}$-isomorphic to $\PP^2$ blown up at the base locus of $\Gamma$. }
   \end{dfn}
 The choice of the pencil $\Gamma$ is non-canonical.

 Suppose $\Gamma$ induces $\cc{E}$ and $p_1,...,p_r$ are the distinct base points of $\Gamma$. The Picard group of $\cc{E}$ is generated by the strict transform of a cubic in $\Gamma$, the exceptional curves above each $p_i$ and some of the $(-2)$-curves obtained in the process of blowing up the $p_i$'s in case the points have multiplicity strictly larger than one as base points or are in a non-general position (for example, three on a line or six on a conic).  As the exceptional curves are the sections of the elliptic fibration, the above gives us the following information about the (geometric) Mordell-Weil group of $\cc{E}$.
 
 \begin{lem}\label{mwrk} Let $\Gamma$ be a pencil of cubics in $\PP^2$ and $\cc{E}$ the elliptic surface induced by $\Gamma$. Let $s$ be the number of distinct base points of $\Gamma$ and $\bar{r}$ be the geometric Mordell-Weil rank of $\cc{E}$, i.e., the rank of the Mordell-Weil group over the algebraic closure $\bar{k}$. Then \[\bar{r}\leq s-1.\]
 \end{lem}

 Over a number field, a rational elliptic surface may have a minimal model other than $\mathbb{P}^2$. Hence we cannot assure the existence of a
 $k$-birational morphism between $\cc{E}$ and $\PP^2$. We treat this situation in Subsection \ref{minmodarith}.
 
\subsection{N{\'e}ron-Tate height in elliptic surfaces}

 In \cite{S2} Shioda developed the theory of Mordell-Weil
 lattices. He remarked that the N{\'e}ron-Severi and the Mordell-Weil groups modulo torsion
 have a lattice structure endowed with a pairing given essentially by
 the intersection pairing on the surface. In the Mordell-Weil group it coincides
 with the N{\'e}ron-Tate height. 

 To define this pairing we need to introduce some notation:

 Let $\Theta_v$ be a fibre with $m_v$ components denoted by
 $\Theta_{v,i}$. Let $\Theta_{v,0}$ be the zero component (the one
 that intersects the zero section) and $A_v=
 ((\Theta_{v,i} .\Theta_{v,j}))_{1\leq i,j\leq m_v-1}$ a (negative
 definite) matrix. The pairing is given by the following formula:

 \[ \langle P,Q \rangle = \chi +(P.O)+(Q.O)-(P.Q)-\sum_{v \in R} contr_v(P,Q),\]
where $\chi$ is the Euler characteristic of the surface, $R$ is the set of reducible fibres, $(P.Q)$ the
 intersection of the sections given by $P$ and $Q$, and $contr_v(P,Q)$
 gives the local contribution at $v$ according to the intersection of
 $P$ and $Q$ with the fibre $\Theta_v$: if $P$ intersects
 $\Theta_{v,i}$ and $Q$ intersects $\Theta_{v,j}$ then $contr_v(P,Q)=
 -(A_v^{-1})_{i,j}$ if $i,j \geq 1$ and zero if one of the sections
 cuts the zero component. For a table of all possible values of $contr_v$ according to the
 fibre type of $\Theta_v$ see \cite{S2}.
 
 We will only use the fact that the N\'{e}ron-Tate height of a section can be computed in terms of its numerical class.
     
\section{Arithmetic preliminaries}\label{arith}
\subsection{Minimal models over perfect fields}\label{minmodarith}

 The theory in this subsection was developed by Enriques, Manin \cite{Ma1},
 \cite{Ma2} and Iskoviskikh \cite{I}. We state the main results. For proofs
 we invite the reader to look at the bibliography cited above.
 
 \begin{thm} Let $X$ be a smooth minimal rational surface defined over a perfect field k and let $\mathrm{Pic}(X)$ denote its Picard group over $k$. Then $X$ is isomorphic to a surface in one of the following families:
 \begin{itemize}
 \item[I)] A del Pezzo surface with $\mathrm{Pic} (X)\simeq \ZZ$.
 \item[II)] A conic bundle such that $\mathrm{Pic}(X) \simeq \ZZ \oplus \ZZ$.
 \end{itemize} 
 Reciprocally, if X belongs to family I) then it is (automatically) minimal. If X belongs to family II) then it is not minimal if, and only if, $d=3,5,6$ or $d=8$ and $X$ is isomorphic to the ruled surface $\mathbb{F}_1$. There are no minimal surfaces with $d=7$.
 
 Some surfaces endowed with a conic fibration are at the same time del Pezzo surfaces, namely, if $d=3,5,6$ or $d=1,2,4$ and X has two distinct conic fibrations, or $d=8$ and $\bar{X}= X \times_k \bar{k} \simeq \PP^1 \times\PP^1$ or $\PP^2$ blown up in one point.
  \end{thm}
 
 \begin{dfn} \emph{We say that a surface $X$ is $k$-birationally
   trivial or $k$-rational if there is a birational map $\PP^2 \dashrightarrow X$ defined over $k$.}
 \end{dfn}

\begin{thm} Every minimal rational surface with $d\leq 4$ is $k$-birationally non-trivial.
\end{thm}

\begin{thm}\label{deg3}
Every rational surface X of degree at least five and such that $X(k)\neq \emptyset$ is $k$-rational and every rational surface $X$ with $d \geq 3$ and  $X(k) \neq \emptyset $ is $k$-unirational.
\end{thm}

\begin{rem}\emph{A priori most $k$-minimal surfaces with $d\geq 1$ and $X(k)\neq \emptyset$ can be a $k$-minimal model of a rational elliptic surface. The condition $X(k)\neq \emptyset$ comes from the zero section that is defined over $k$ and is contracted to a $k$-rational point. We will exclude the conic bundles such that $\bar{X}= X \times_k \bar{k}$ is isomorphic to $\PP (\cc{O}_{\PP^1}\oplus \cc{O}_{\PP^1}(n))$ with $n \geq 3$ since these surfaces contain a curve with self intersection $-n$ and a rational elliptic surface contains no curves with self intersection $-n$ for $n\geq 3$.} 
\end{rem}

\begin{rem}
 \emph{If a surface $X$ of degree $d=\omega^2_X$ is a $k$-minimal model of a rational elliptic surface $\cc{E}$, then $\cc{E}$ is isomorphic to the blow up of $X$ in $d$ points which form a Galois invariant set.}
\end{rem}

We finish this subsection giving a sufficient condition for $\bar{k}$-rational elliptic surfaces to be $k$-unirational.

\begin{lem}\label{unirational}
Let $\cc{E}$ be a $\bar{k}$-rational elliptic surface defined over a perfect field $k$. If $\mathrm{rank}(\mathrm{Pic}(\cc{E}|k)) \geq 5$ then $\cc{E}$ is $k$-unirational.
\end{lem}

\noindent\textbf{Proof:} Since the Picard group of $k$-minimal rational surfaces has rank one or two, any $k$-minimal model of $\cc{E}$ will be obtained after contracting at least three $(-1)$-curves. Hence if $X$ is a $k$-minimal model of $\cc{E}$ then $K_X^2\geq 3$. Moreover, $X(k) \neq \emptyset$ since the image of the zero section will provide at least one $k$-rational point. It follows then from Theorem \ref{deg3} that X is $k$-unirational. \qed

\subsection{Kummer theory}\label{kummer}

 Let $K$ be a number field or a function field and $A$ an abelian
variety. Let $P\in A(K)$ be a point of infinite order. We say that $P$ is \emph{indivisible by $n$} if for all $Q\in
A(K)$ such that there exists a divisor $d$ of $n$ such that $[d]Q=P$ we have $d= \pm 1$.
\begin{rem}
\emph{In \cite{M}, Hindry defines a point as being \emph{indivisible} if it is indivisible by all natural numbers $m >1$. One can replace this definition by the one introduced above, i.e., $P$ indivisible by $m$ in the hypothesis of \cite[Lemme 14]{M}. In fact, let $P$ be a point indivisible by $m$. We can write $P=[l]P_1$ with $P_1$ indivisible as in \cite{M} and $(l,m)=1$. There exist $u, v \in \mathbb{N}$ such that $ul+vm=1$ which allows us to write \[P_1=[u]P+ [vm]P_1,\] and thus \[ K(\frac{1}{m}P) \subseteq K(\frac{1}{m}P_1)=K(\frac{u}{m}P) \subseteq K(\frac{1}{m}P).\] Hence $K(\frac{1}{m}P_1)= K (\frac{1}{m}P)$.
In the rest of this subsection we state some of the results that can be found in \cite{M} about the degree of the extension $K(\frac{1}{m}P)$ taking this remark into account, i.e., replacing the hypothesis $P$ indivisible, by the hypothesis $P$ indivisible by $m$.}
\end{rem}
 
 Let $n\in \NN$ and $P$ be a point indivisible by $m$ where $m$ divides $n$. Denote by $A_n$ the set of $n$-torsion points in $A(\bar{K})$ and
 $\frac{1}{m}P$ a point $Q \in A(\bar{k})$ such that $[m]Q=P$. The Galois group $G_{n,K,\frac{1}{m}P}$ of
 the extension \[ K(A_n, \frac{1}{m}P)|K(A_n)\] can be viewed as a subgroup of $A_m$. Kummer theory for abelian
 varieties tells us that if $K$ is a number field the group $G_{n,K,\frac{1}{m}P}$ is actually almost the whole group $A_m$, i.e., its image inside $A_m$ by an injective map is of
 finite index which is bounded by a constant independent of $P$ and $m$.

 This theory was studied in full generality by Ribet, Bertrand,
 Bashmakov and others. We restrict ourselves to elliptic
 curves. For us $k$ is a number field and $K=k(B)$ is the function
 field of a projective curve $B$.

 We now state a typical result of Kummer theory for elliptic
 curves. This is valid in a more general context (see \cite{Rib}). 

\begin{thm}\label{nbf}
Let $E$ be an elliptic curve defined over a number field $k$, $n \in \NN$ and $P
\in E(k)$ a point indivisible by $m$ where $m$ divides $n$. There exists a positive constant $f_0=f_0(E,k)$ such that\[\ | G_{n,k,\frac{1}{m}P}| \geq f_0.m. \]
\end{thm} 

\noindent \textbf{Proof:} See \cite[Appendix 2]{M} for the case when $P$ is an indivisible point.
\indent\par
\indent\par
 Let $\pi: \cc{E} \rightarrow B$ be an elliptic surface and $P$ a point indivisible by $n$ in the generic fibre. Since Galois groups become smaller after specialization, after applying the previous theorem to a fibre $E= \pi^{-1}(t)$ defined over the number field $k$ we have \[|G_{n, k(B),\frac{1}{m}P}| \geq | G_{n,k,\frac{1}{m}P}| ,\] which gives us a theorem as Theorem \ref{nbf} for elliptic curves over function fields.
 
 \begin{thm}\label{ff}
Let $E$ be an elliptic curve defined over a function field $K$, $n\in \NN$ and $P
\in E(K)$ a point indivisible by $m$ where $m$ divides $n$. There exists a positive constant $f_0=f_0(E,K)$ such that\[\ | G_{n,K,\frac{1}{m}P}| \geq f_0.m .\]
\end{thm}

\begin{rem}
\emph{ If the surface $\cc{E}$ above is not isotrivial then the endomorphism ring of the generic fibre $E$ is $\ZZ$ and the result above is stronger, namely, the group $G_{n,K,\frac{1}{m}P}$ is almost all the set of $m$-torsion points $\cc{E}_m$ i.e., its image inside the $m$-torsion subgroup under an injective map is of finite index, bounded by a constant independent of $P$ and $m$ (c.f \cite[Proposition 1]{M}).} 
\end{rem}

We finish this subsection with a lemma about torsion points on elliptic curves that will be used in the next section during the proof of Proposition \ref{prpimportante}. See \cite{Serre} for a proof and more results on torsion points on elliptic curves.

\begin{lem}\label{serre}
Let $E$ be an elliptic curve defined over a field $k$, which is either a number field or a function field, and $P \in E(\bar{k})[m] \setminus E(k)$, a point of order $m$. Then there exist an $\alpha>0$ and a constant $c_{E, k}$ independent of the point $P$ such that \[[k(P):k]\geq c_E.m^{\alpha}.\]
\end{lem}
\section{Proof of Theorem \ref{th2}}
Let $\pi:\cc{E}\rightarrow B$ be as in Section 2. In the first subsection we show that given a non-constant pencil of curves in $\cc{E}$ that are not contained in a fibre of $\pi: \cc{E}\rightarrow B$, then all but finitely many curves in it yield, after base change (see Subsection \ref{bc}), a \textit{new section} that is independent of the \textit{old sections}. The proof of Theorem $\ref{th2}$ then depends only on a construction of a family of irreducible curves defined over $k$ with infinitely many $k$-rational points, i.e., $\mathbb{P}^1$ or genus-one curves with positive Mordell-Weil rank. This construction will be given in the second subsection.

\subsection{A key proposition}

 First we state the most useful criteria for us to determine when a \textit{new section} is independent of the \textit{old} ones. 
 
 \begin{lem}  An irreducible curve $C \subset \cc{E}$ that is not a component of a fibre induces a new section on $\cc{E}\times_B C$ independent of the \textit{old} ones if and only if for every section $C_0 \subset \cc{E}$ and every $n \in \NN^*$, the curve $C$ is not a component of $[n]^{-1}(C_0)$.   
 \end{lem}

 Let $\cc{C}$ be a family of curves in a projective surface $X$. We call $\cc{C}$ a \emph{numerical family} if all its members belong to the same numerical class in the N\'eron-Severi group. We prove that in a numerical family it is enough to check the conditions of the lemma for a bounded $n$ and a finite number of sections. Thus if the family is infinite, all but finitely many members induce new independent sections after base change.
 
 \begin{pro}\label{prpimportante}
 Let $\cc{E}\rightarrow B$ be an elliptic surface defined over a number field $k$. Let $\cc{C}$ be a numerical family of curves inside $\cc{E}$. There exist an $n_0(\cc{C}) \in \NN$ and a finite subset $\Sigma_0(\cc{C})\subset
 Sec(\cc{E})$ such that for $C$ in the family $\cc{C}$, the new section induced by C
 is linearly dependent of the old ones if, and only if, $[n]C \in \Sigma_0$ for some $n\leq n_0$.  
 \end{pro}

\noindent \textbf{Proof:}  Suppose $[n]C= C_0$ for some section $C_0$. We may assume such $n$ to be minimal i.e., there does not exist $n'<n$ such that the curve $[n']C$ is a section. 
The proof is divided in two parts:
\begin{itemize}
\item[1)] Bounding $n$ from above using Kummer theory (see Section \ref{kummer}).
\item[2)] For a fixed $n$ such that $C_0=[n]C$, showing that the set of sections in the same numerical class is finite by N\'eron-Tate height theory.
\end{itemize}

\textbf{1) Bounding n:} We define the degree of a curve in $\cc{E}$ by its intersection with a fibre:
\[ \mathrm{deg}(C)= (C.F).\] 

 If $C_0$ is a section we have $\mathrm{deg}(C_0)=(C_0.F)=1$. The degree of $C$, which will be denoted by $h$, is fixed within the family, since all curves belong to the same numerical class.
 
  The map $[n]$ is not a morphism defined on the whole surface, but on an open set $U \subseteq \cc{E}$ which excludes the singular points in the fibres. Since sections do not intersect the fibres in singular points they are contained in $U$. This allows us to write:
  \[\mathrm{deg}([n]^{-1}C_0)= (([n]^{-1}C_0).F)= n^2(C_0.F)= n^2.\]

 Thus, $\lim_{n \to \infty}\deg [n]^{-1}(C_0) = \infty$.
 
 Denote by $K$ the field $k(B)$, by $E$ the generic fibre of $\cc{E}$ and by $P_0$ the point in $E(K)$ corresponding to the section $C_0$.  Let $P\in E(\bar{K})$ be such that \[ [n] P=P_0 \] where $n$ is minimal with respect to the expression above. We now show that $n$ is bounded by a constant $n_0$ that depends only on $E$, $K$ and $h$.
 
  Note first that if $P$ is a torsion point, we know by Lemma \ref{serre} that its order $m$ is bounded by $ c_E.m^{\alpha}\leq [K(P):K]= h$.
  
Now suppose $P$ is of infinite order. Let $m$ be the smallest positive integer such that there exist $P_1 \in E(K)$ and $T$ a torsion point satisfying \[ [m] P=P_1 + T.  \]
 We claim that $P_1$ is indivisible by $m$. If $l$ is a divisor of $m$ such that $[l]Q = P_1$ with $Q \in E(K)$, then \[[l][\frac{m}{l}]P=[l]Q +T\] and hence there exists a torsion point $T_1$ such that \[[\frac{m}{l}]P= Q+T_1.\] By the minimality of $m$ with respect to the equation above we must have $l=\pm 1$. 
 
  Let $m'=mm_1$ where $m_1$ is the order of the torsion point $T$. Let $T'$ be such that $[m]T' = T$ and put 
$P' = P + T'$ ; then $[m]P' = P_1$ and $T' \in E_{m'}$ . Applying Theorem \ref{ff} we obtain
  \begin{equation}\label{kum}  
   [K(P',E_{m'}): K(E_{m'})] \geq m.f_1.
  \end{equation}
  Now, note that $K(P, E_{m'} ) = K(P' , E_{m'} )$ hence \[ h=[K(P ) : K] 
\geq [K(P, E_{m'} ): K(E_{m'} )] = [K(P' , E_{m'} ): K(E_{m'} )] \geq f_1 .m,\] 
and thus m is bounded. Since $T$ is deÞned over $K(P )$, its order $m_1$ is also bounded in terms of $h$, and thus $[mm_1] P = 
m_1 P_1\in E(K)$ with $n\leq mm_1$ bounded as announced. 
  
  \textbf{2) The numerical class of a section:}
   Fix $C_1,...,C_r$ generators of the Mordell-Weil group.
   
   For a fixed $n$, the intersection multiplicity $(([n]C).C_i)$ is also fixed, say equal to $n_i$, and it depends only on the numerical class of $C$. The same holds for the intersection of $C$ with the zero section, say equal to $m_0$ and for the intersection with the fibre components $\Theta_v$, say equal to $ l_{v,j_v}$. 
    The N\'{e}ron-Tate height in an elliptic surface is uniquely determined by the intersection numbers above. The set
  \[\Sigma_0= \{ C_0 \text{ section } | (C_0. C_i)= n_i; i= 1,..., r, \,\
(C_0. O)= m_0 \, \mathrm{and} \, (C_0. \Theta_{v,j_v})= l_{v,j_v}\}\] 
is finite since it is a set of points with bounded N\'eron-Tate height, thus there are only finitely many possible sections $C_0$ such that $C \subset [n]^{-1}C_0$.
\qed

\begin{cor}\label{principal} 
Let $\cc{E}$ be an elliptic surface and $\cc{L}$ a non-constant numerical family of curves on $\cc{E}$ whose members are not contained in the fibres of $\cc{E}$, then for almost all member $C$ of the pencil, the new section induced by $C$ is independent of the old sections.
\end{cor}

\subsection{Proof of Theorem \ref{th2}}

Let $\psi:\mathbb{P}^2 \dashrightarrow \cc{E}$ be a $k$-unirational map. Let $\cc{L}$ be given by the set of lines in $\mathbb{P}^2$. Then $\cc{L}^{\psi}= \{\psi(L); L \in \cc{L}\}$ is an infinite family of curves in $\cc{E}$ defined over $k$ whose general member is integral and of geometric genus zero. These curves cannot be all contained in a fibre of $\cc{E}$ as the family is infinite, with irreducible members, and there is only a finite number of reducible fibres for the elliptic fibration in $\cc{E}$.

The theorem follows from an application of Corollary \ref{principal} to the family $\cc{L}^{\psi}$. \qed

\section{Proof of Theorems \ref{th1} and \ref{th3}}
 
  To prove Theorems \ref{th1} and \ref{th3} we need to produce two families of curves defined over $k$. These families must, not only have infinitely many $k$-rational points, but also be such that the fibred product of two curves in different families is irreducible. The first subsection is devoted to the construction of such families, first done in the context of Theorem \ref{th1}, i.e., over $\mathbb{P}^2$; then in the settings of Theorem \ref{th3}, i.e., over other $k$-minimal surfaces. The proof of irreducibility of the generic member of the constructed families is given in the second subsection and is followed by the verification that such curves contain indeed infinitely many $k$-rational points. Finally, all is assembled in Subsection \ref{prf} to conclude the proofs of Theorems \ref{th1} and \ref{th3}.   

\subsection{Construction of linear pencils of rational curves}\label{constructions}

The results presented in this subsection are technical and may be skipped by the reader willing to accept the existence of two linear pencils of conics, i.e., $\mathbb{P}^1$'s intersecting the fibres with multiplicity two in the surfaces satisfying the hypothesis of Theorem \ref{th1}. Here we provide ``case-by-case", depending on the configuration of blown up points, constructions of linear pencils of curves on $\cc{E}$ to which we apply Corollary \ref{principal}. We construct two linear pencils of rational curves defined over $k$ in a $k$-minimal model of the rational elliptic surface $\cc{E} \rightarrow B$. Since the base-changed surface fibres over the fibred product of these two curves over $B$, we fabricate those curves in a way that their fibred product has genus at most one. We state below sufficient conditions for this.

\begin{lem}\label{deg2} Let $C_1$ and $C_2$ be two smooth projective rational curves given with two distinct morphisms $\varphi_i: C_i \rightarrow B$ of degree 2 to a genus zero curve $B$. Then $g(C_1\times_B C_2)\leq 1$.
\end{lem}
\noindent \textbf{Proof}: It is a simple application of Hurwitz' Formula. \qed

 We can proceed to the constructions. They depend on the degree of the $k$-minimal model considered, as well as on the configuration of the blown up points under the action of the absolute Galois group $\mathrm{Gal}(\bar{k}|k)$. We start with the simplest case, namely, when $\cc{E}$ has a minimal model isomorphic to $\mathbb{P}^2$.
\\
\\
\textbf{a) A minimal model $k$-isomorphic to $\PP^2$}

 Let $p_1,..,p_9$ be the nine, not necessarily distinct, points in the
blow up locus. Since the zero section is defined over $k$, at least one
of the points above is $k$-rational, say $p_1$. If $C$ is a curve of degree $d$ in $\mathbb{P}^2$ with multiplicity $m_i$ through $p_i$ then its genus satisfies $g(C) \leq \frac{(d-1)(d-2)}{2}- \sum_i \frac{m_i(m_i-1)}{2}$. Denote by $C'$ its strict transform under the blow up $\cc{E} \rightarrow \mathbb{P}^2$, then the degree of the map given by the restriction of $\pi$ to $C'$ is given by $\mathrm{deg}(C'\rightarrow B)= (F. C')= 3d- \sum_i m_i$ where $F$ is a fibre of the elliptic fibration $\pi$. 

Let $\cc{L}_1$ be the pencil of lines in $\PP^2$ through $p_1$. Let $\cc{L}'_1$ be the
pencil of curves in $\cc{E}$ given by the strict transforms of the curves in $\cc{L}_1$ and $C'$ a curve in
$\cc{L}'_1$-- from now on the superscript $'$ will denote the pencil
or curve in the elliptic surface given by the strict transform of that in the minimal model. Then \[\mathrm{deg}(C' \rightarrow B)= 3.1-1=2.\]
We now construct a second pencil of rational curves $\cc{L}_2$ such that the curves in the pencil of strict transforms induced in $\cc{E}$ satisfy Lemma \ref{deg2}.

Since the blow up is defined over $k$, the set formed by the other points is invariant
under the action of $\mathrm{Gal}(\bar{k}|k)$. The construction depends on the size of
the smallest orbit different from $p_1$ whose points are, together with $p_1$, in general position, i.e., no three are collinear, no six lie in a conic and there is no cubic through eight of the points singular in one of them.
\\
 \begin{itemize}
\item[i)] One other $k$-rational point $p_2$.

In this case we can construct $\cc{L}_2$ in a similar way as we did for $\cc{L}_1$: take \[ \cc{L}_2= \{l \text{ a line in } \PP^2 \text{ through } p_2\} .\]
Any curve in $\cc{L}'_2$ together with any curve in $\cc{L}'_1$ satisfy Lemma \ref{deg2}.
\\
\item[ii)] Two conjugate (under $\mathrm{Gal}(\bar{k}|k)$) points $p_2,p_3$.

Let $\Lambda$ be a pencil of cubics in $\mathbb{P}^2$ inducing $\cc{E}$ such that $p_1,p_2,p_3$ are base points of it. Since we suppose that points are in general position, $p_1, p_2$ and $p_3$ are not collinear. Let us first suppose that there are no other base points and thus that the multiplicities $(m_1,m_2,m_3)$ of $p_1,p_2,p_3$ as base points of $\Lambda$ are $(1,4,4), (3,3,3), (5,2,2)$ or $(7,1,1)$. In the first case, every cubic in $\Lambda$ shares the same tangent lines, say $l_2$, through $p_2$ as well as the same tangent line, say $l_3$, through $p_3$. We consider $\cc{L}_2$, the set of conics through $p_2,p_3$ with tangents $l_i$ through $p_i$, for $i=2,3$. Let $C$ be a conic in $\cc{L}_2$, then C intersects the cubics of $\Lambda$ in $p_2$ and $p_3$ with intersection multiplicity two and in two other points. Hence, the morphism from the strict transform $\varphi_{C'}:C' \rightarrow B$, given by the restriction of the fibration to $C'$, has degree two. 
In the remaining three cases all cubics of $\Lambda$ share a tangent line, say $l_1$, through $p_1$. We take $\cc{L}_2$ to be the set of conics through $p_1,p_2,p_3$ with prescribed tangent $l_1$ through $p_1$. It is a linear pencil of conics, since the space of conics in $\mathbb{P}^2$ has dimension 5.  A conic $C$ in $\cc{L}_2$ also intersects the cubics in $\Lambda$ at the point $p_1$ twice, at the points $p_2,p_3$ and at two other points. Thus if $C$ is a conic in $\cc{L}_2$ then the morphism $\varphi_{C'}$ from $C'$ to $B$ has also degree two.
As in i) any curve in $\cc{L}'_2$ together with any curve in $\cc{L}'_1$ satisfy the hypothesis of Lemma \ref{deg2}.

Now suppose there are other base points. If there is another orbit with two or four points one may apply construction iv) below independently of the configuration of these extra points. If there are three, five or six other conjugate base points then apply construction iii), v), vi), respectively, below. 
 \\
\item[iii)] Three conjugate points $p_2,p_3,p_4$.

The construction here is simpler. The pencil of conics $\cc{L}_2$ through $p_1,p_2,p_3$ and $p_4$ is such that for every curve $C \in \cc{L}_2$ the morphism $\varphi_C'$ has degree two.
\\
\item[iv)] Four conjugate points $p_2,p_3,p_4,p_5$.

As in the previous case, but now we let $\cc{L}_2$ be the pencil of conics through $p_2,p_3,p_4$ and $p_5$.
\\
\item[v)] Five conjugate points $p_2,p_3,p_4,p_5,p_6$.

 Take $\cc{L}_2$ to be the pencil of cubics through $p_1,...,p_6$ with a singularity at $p_1$. The degree of the morphism $\varphi_C$ for $C \in \cc{L}_2$ is two $(9-5-2=2)$. So that together with a curve in $\cc{L}'_1$ the curve $C$ satisfies the hypothesis of Lemma \ref{deg2}.
\\
\item[vi)] Six conjugate points $p_2,...,p_7$.

 Consider $\cc{L}_2$ to be the pencil of quintics singular at each $p_2,...,p_7$ that passes through $p_1$ as well. Since $\mathrm{dim}(H^0(\PP^2, \cc{O}(5)))=21$ and we have imposed 19 conditions, $\cc{L}_2$ forms at least a linear pencil. The curves in it are rational since $ g\leq 6-6=0.$
 
  The degree of $\varphi_C'$ for $C \in \cc{L}_2$ is equal to $5.3-6.2-1=2$.
\\
\item[vii)] Seven conjugate points $p_2,...,p_8$.

Consider $\cc{L}_2$ the space of quartics through the seven points $p_2,...,p_8$ with multiplicity at least three at $p_1$. That gives at most 13 conditions in a space of dimension 15. So $\cc{L}_2$ is at least a linear pencil. Its curves are rational and the degree of the induced morphism to $B$ is equal to $4.3-7-3=2$.
\\
\item[viii)] Eight conjugate points $p_2,...,p_9$.

We consider highly singular curves. Take $\cc{L}_2$ to be the set of curves of degree 17 such that
\begin{itemize}
\item The eight points $p_2,...,p_9$ are singular with multiplicity at least six.
\item It passes through $p_1$.
\end{itemize}
This gives us at most 169 conditions in a 171-dimensional space and thus at least a linear pencil. The degree of the morphism is also two and the curves have genus zero as in the previous cases.
\end{itemize}
\textbf{}

\textbf{Irreducibility of the curves constructed above:}  Since the points considered are in general position the curves constructed above are irreducible. This is trivially verified in cases i),..., iv). In cases v),...,viii) one can easily check that if the pencil is generically reducible then either its base points are in non-general position i.e., there are three collinear points, six lie on a conic or eight on a cubic singular at one of them, or the Galois orbits break into smaller orbits.
\indent\par
\indent\par
We now focus on the other possible minimal models over $k$. As observed after the statement of Theorem \ref{th3}, we may suppose $X$ is a del Pezzo surface. Let us first recall some geometric and arithmetic facts about those surfaces.

 Since we are dealing with surfaces whose Picard group over $k$ is
 small (isomorphic to $\ZZ. \omega_X$), the most natural place to look for curves defined over $k$ is $H^0(X, \omega^{-n}_X)$. We now recall the dimension of these spaces and the genus of the curves in them.
 \begin{lem}\label{dimdel} Let $X$ be a $k$-minimal del Pezzo surface of degree $ 1\leq (\omega_X. \omega_X)=d \leq 8$ defined over a number field $k$. Given points $p_1,...,p_j$ in $X(\bar{k})$ and non-negative integers $n_1,...,n_j$, let $\cc{L}=\{ s \in H^0(X, \omega^{-n}_X); m_{p_i}\geq n_i\}$ where $m_{p_i}$ denotes the multiplicity at the point $p_i$ of the curve given by the divisor of zeros of the section $s$. The following hold:
 \begin{itemize}
\item[i)] $\mathrm{dim} (\cc{L})\geq \frac{d(n^2+n)}{2}+1-\sum_i \frac{n_i^2+n_i}{2}$.
\item[ii)] If $L \in \cc{L}$ then $g(L)\leq \frac{d(n^2-n)}{2}+1-\sum \frac{n_i^2-n_i}{2}$.
\item[iii)] If $\pi: \cc{E}\rightarrow B$ is an elliptic surface obtained by blowing up $p_1,...,p_j$ and $L'$ is the strict transform of $L$ in $\cc{E}$ then $\mathrm{deg} (\pi|_{L'}: L' \rightarrow B)=nd-\sum m_{p_i}$.
\end{itemize}
 \end{lem}
\noindent \textbf{Proof:} See \cite[Chapter III Lemma 3.2.2]{Kollar}. 
 \indent\par
\indent\par
We can now proceed to the construction of linear pencils on the
$k$-minimal models. We recall that at least one point in the blow up
locus of $f:\cc{E}\rightarrow  X$ is $k$-rational, the one that comes
from the contraction of the zero section. We suppose that $\cc{E}$ has no reducible fibres. It follows that all the
points on the locus of $f$ are distinct since the blow up of infinitely near points gives rise to $(-2)$-curves, and these are always components of reducible fibres. For (i) of Theorem \ref{th3} it is clearly
sufficient to do the construction in the case where there are two
orbits by the action of the Galois group: the one of the $k$-rational point and
another one with the other $d-1$ points. As in the case where $\mathbb{P}^2$ was a $k$-minimal model, we
will look for curves satisfying the hypothesis of Lemma \ref{deg2}.
\\
\\
\textbf{b) A minimal model isomorphic to a del Pezzo surface of degree eight:}

 Let $p_2,...,p_9$ be the points on the blow up locus of $f:\cc{E} \rightarrow X$. Let $p_2$ be the $k$-rational point. We consider the following pencils of curves in $X$:

\[ \cc{L}_1= \{ s \in H^0(X,\omega^{-8}_X); m_{p_2}(s) \geq 13, m_{p_i}(s)\geq 7, i=3,...,9\} \]
and
\[ \cc{L}_2= \{s \in H^0(X, \omega^{-22}); m_{p_2}(s)\geq 13, m_{p_i}(s)\geq 23, i=3,...,9\}. \]

 If $C'_1 \in \cc{L}'_1$ then by Lemma \ref{dimdel} \[ g(C'_1)\leq 8 \frac{(64-8)}{2}+1-\frac{(169-13)}{2}-7\frac{(49-7)}{2}=0 \text{ and } \mathrm{deg}(\varphi_{C_1})=64-13-49=2.\]

 If $C'_2 \in \cc{L}'_2$ then \[g(C'_2)\leq 8\frac{(22^2-22)}{2}+1-\frac{(169-13)}{2}-7\frac{(23^2-23)}{2}=0 \text{ and } \mathrm{deg}(\varphi_{C_2})=22.8-13-7.23=2.\]

Thus $C'_1$ and $C'_2$ satisfy the hypothesis of Lemma \ref{deg2}. Since, by Lemma \ref{dimdel} (i), they belong to a linear pencil of curves, by Corollary \ref{principal} they can be chosen in a way such that the new sections induced by them in the base-changed surface are independent of the \textit{old sections} and of each other.  
 
 Since there are no minimal rational surfaces of degree seven we now pass to surfaces of degree six.
\\
\\
\textbf{c) A minimal model isomorphic to a del Pezzo surface of degree six:}

 Let $p_4,...,p_9$ be the points on the blow up locus of $f$. We consider the possible orbits under the action of the absolute Galois group. We denote the cases by $(n_1,...,n_r)$ where $r$ is the number of distinct orbits and $n_i$ is the multiplicity of the points in the same orbit.
 \begin{itemize}
 \item[i)] If the points lie in a $(1,2,3)-$configuration then the blow up of the two points in the same orbit produces a surface of degree four to which we apply the constructions in e). 
\item[ii)]  If $(1,1,n_3,n_4)$ let $p_4$ and $p_5$ be the $k$-rational points. The blow up of $p_4$ produces a surface of degree five to which we can apply the constructions in d). 
 \end{itemize}
\noindent 
 \textbf{d) A minimal model isomorphic to a del Pezzo surface of degree five:}

 Here we consider the following pencils:
 
  \[\cc{L}_1= \{s \in H^0(X, \omega_X^{-2}); m_{p_5}(s)
 \geq 4, m_{p_i}(s)\geq 1, i=6,...,9 \}\] and \[ \cc{L}_2= \{s \in H^0(X, \omega_X^{-10}); m_{p_5}(s) \geq 4, m_{p_i}(s) \geq 11, i=6,...,9\}.\]

\[\mathrm{dim}(\cc{L}_1)\geq16-10-4=2, \,\, g(C_1)=6-6=0, \,\, \mathrm{deg}(f:C'_1 ,\rightarrow D)=10-4-4=2,\]
\[\mathrm{dim}(\cc{L}_2)\geq 5(110)/2+1-10-4(66)=2.\]

\noindent
\textbf{e) A minimal model isomorphic to a del Pezzo surface of degree four:}

 The pencils that we consider to prove Theorem \ref{th1} are:
 \[\cc{L}_1= \{s \in H^0(X, \omega^{-1}_X); m_{p_1}(s)\geq 2\}, \]

 and

\[\cc{L}_2= \{s
 \in H^0(X, \omega^{-7}_X); m_{p_1}(s) \geq 2 , m_{p_i}(s) \geq
 8, i=2,3,4\}.\]

\subsection{Irreducibility of the curves}

We prove that if $\cc{E}$ satisfies the hypothesis of Theorem \ref{th3} the curves constructed in the previous subsection are irreducible. 
First, we show that a series of Cremona transformations, i.e., the blow up of three distinct and non-collinear points followed by the contraction of the three lines through them, reduces each pencil of curves produced above to a pencil of lines in $\mathbb{P}^2$ passing through a point. Cremona transformations are not in general automorphisms of the
surface, but they are automorphisms of the Picard group. In particular a class is represented by a connected curve if and only if the transformed
one under a Cremona transformation is represented by an connected curve. We then use the fact that $\cc{E}$ has no reducible fibres to show that the curves constructed are irreducible (See \cite{testa} for more on irreducibility of spaces of curves on del Pezzo surfaces).

\begin{lem}
  Let $\cc{E}$ be a rational elliptic surface with no reducible fibres, $X$ a $k$-minimal model of $\cc{E}$ as in Theorem \ref{th3} and let $\cc{L}$ be one of the pencils of curves constructed in the previous subsection. Then the generic member of $\cc{L}$ is an irreducible curve with geometric genus zero.
\end{lem}
 
\noindent \textbf{Proof:} Let $d$ be the degree of $X$. Let $f: \cc{E} \rightarrow X$ be the map corresponding to the blow up of a $\mathrm{Gal}(\bar{k}|k)$-invariant set of distinct points $P_1,...,P_d \in X(\bar{k})$. Since $\cc{E}$ has no reducible fibres, the exceptional curves above  $P_i,$ for $i=1,...,d$, are all independent in the Mordell-Weil group of $\cc{E}$. Therefore, they provide a subset of a set of generators of the Picard group of $\cc{E}$.
We can fix a basis for the geometric Picard group of $\cc{E}$ to be $\{L_0,...,L_9\}$ where $L_0$ is the total transform of a line $l$ in $\mathbb{P}^2$, $L_1,...,L_d$ are the exceptional curves above $P_1,...,P_d$ and $L_{d+1},...,L_9$ are also exceptional curves in $\cc{E} \times_k{\bar{k}}$.

Let $g: \cc{E} \rightarrow \mathbb{P}^2$ be a blow up presentation, defined over $\bar{k}$, factoring through $f$. We represent a curve $C$ in $\cc{E}$ by its \emph{numerical type}, i.e., by the list of coordinates of its divisor class in the basis given by $\{L_0,...,L_9\}$ of the Picard group:
\[(d,m_1,...,m_9), \]
where $d$ is the degree of the image of $C$ in $\mathbb{P}^2$ with respect to $g: \cc{E} \rightarrow \mathbb{P}^2$,  $m_i=m_{P_i}(C)$, the multiplicity of the curve $C$ at the point $P_i \,, i=1,...,9$. 
 
  Let $C$ be a curve in $X$, then the strict transform of $C$ through $f$ is a curve in $\cc{E}$ given by \[f^{-1}(C)- \sum_{i=1,...,d} m_{p_i}(C)L_i.\]  We define the \emph{numerical type} of $C$ as the numerical type of its strict transform in $\cc{E}$. 

For example if $X= \mathbb{P}^2$, lines in $\mathbb{P}^2$ through the point $p_1$ are represented by the class $(1,1,0,...,0)$.  If $X$ is a del Pezzo surface of degree five, curves in $\cc{L}'_1$ where  \[\cc{L}_1= \{s \in H^0(X, \omega_X^{-2}); m_{P_5}(s)
 \geq 4, m_{P_i}(s)\geq 1, i=6,...,9 \} \]
 are represented by $(6,2,2,2,2,4,1,1,1,1)$.
 
  Since Cremona transformations are given by composing birational maps, to show that a curve is connected, it is sufficient to verify that curves constructed in the previous subsection can be (Cremona-)transformed into a curve whose numerical type is one of $(1,1,0,...,0)$,..., $(1,0,...,0,1)$ i.e., into a line through one of the $p_i$. 
  
  After applying successively Cremona transformations, one can check that all curves constructed in the previous subsection are in the same class as a line through a point and are thus connected curves.  
  
  We give below the result of Cremona transformations applied to curves in $\cc{L}'_1$ where $\cc{L}_1= \{s \in H^0(X, \omega_X^{-2}); m_{p_5}(s)
 \geq 4, m_{p_i}(s)\geq 1, i=6,...,9 \}$ and $X$ is a del Pezzo surface of degree five.
 
   Curves in this family are encoded by $(6,2,2,2,2,4,1,1,1,1)$ which after a Cremona transformation become $(4, 0, 0, 2, 2, 2, 1, 1, 1, 1)$, then
$(2, 0, 0, 0, 0, 0, 1, 1, 1, 1)$ which finally, after a last transformation, become
$(1, 0, 0, 0, 0, 0, 0, 0, 0, 1)$.
   
  Now that we have verified that the curves constructed in the previous subsection are connected, we still have to show that the only connected component is irreducible. 
  
    Let $C$ be a curve in one of the families constructed previously. Then $C.F= 2$ and $C$ has geometric genus zero. Thus if $C$ is a reducible curve then it satisfies one of the following:
 \begin{itemize}
 \item[i)] $C= C_1 \cup C_2 \cup (F_1+...+F_m)$ where $C_1$ and $C_2$ are sections and $F_i$ are components of reducible fibres.
 \item[ii)]  $C= D \cup (G_1+...+G_m)$ where $D$ is an irreducible genus zero curve such that $D.F=2$ and $G_j$ are components of reducible fibres. 
 \end{itemize}  
 
 By Proposition \ref{prpimportante}, case i) can only occur for finitely many curves in a numerical family. Thus, we may suppose that the generic members of the families constructed satisfy case ii). As there are only finitely many reducible components of the fibres the curve $D$ in case ii) is such that $\mathrm{dim}|D|\geq 1$. Since $\cc{E}$ has no reducible fibres both pencils $\cc{L}_1$ and $\cc{L}_2$ constructed have irreducible generic members.

  

 \subsection{Infinitely many rational points on the new base}\label{pos}

 In order to prove the corollaries we must show that the base curve of
 the new elliptic surface (the base changed one) has infinitely many
 $k$-rational points. To prove Corollary \ref{cor2} one needs only one
 base change, thus one has to prove that infinitely many among the
 rational curves constructed in the previous section have a
 $k$-rational point. This is assured since the surface $X$ where
 the pencil is constructed, is $k$-unirational \cite[Theorem 3.5.1]{MT} so, in particular, has a Zariski dense
 set of $k$-rational points.

 \begin{lem}
  Suppose $\cc{E}(k)$ is Zariski dense in $\cc{E}$. Let $\{D_t\}_{t \in \mathbb{P}^1}$ be a non-constant pencil of genus zero curves defined over $k$ in $\cc{E}$. Then infinitely many of its curves are $k$-rational. 
 \end{lem}

 For Corollary \ref{cor1} the base curve is in general an elliptic curve. Since our constructions give us families of possible new bases we look at these families as elliptic surfaces and we show that these elliptic surfaces have a non-torsion section.
 
 \begin{thm}\label{posi} Let $\pi: \cc{E}\rightarrow B$ be a rational elliptic
   surface defined over a number field $k$. Let $\cc{L}_1= \{ C_t; t
   \in \PP^1\}$ and $\cc{L}_2= \{ D_u; u \in \PP^1\}$ be two base point free linear
   systems of $k$-rational curves in $\cc{E}$ defined over $k$ with
   $C_t(k), D_u (k) \neq \emptyset$ for infinitely many $t,u$, such that the morphism given by
   the restriction of the elliptic
   fibration to it has degree two. Then for infinitely many $t\in \PP^1(k)$ and infinitely many $u\in \PP^1(k)$ we have $\# (C_t \times_B D_u) (k) = \infty$.
 \end{thm}
 
\noindent  \textbf{Proof:}  Let $C \in \cc{L}_1$ be such that the base-changed elliptic surface $\cc{E}_C= \cc{E} \times_B C$ has generic rank strictly larger than the generic rank of $\cc{E}$. The surface $\cc{E}_C$ admits a second elliptic fibration (since $\cc{L}_2$ is base point free it is not necessary to blow up points to have the fibration), namely, $\cc{E}_C \rightarrow \mathbb{P}^ 1(u)$ where $\mathbb{P}^1(u)$ is the index-set of the pencil of curves $\cc{L}_2$. The fibres of the latter are exactly the curves $D_u \times_B C$. We will show that infinitely many among them have positive rank by showing that this fibration is covered by an elliptic fibration of positive rank. 

The natural morphism $C \rightarrow \mathbb{P}^1(u)$ gives us the surface $\cc{E}_C \times_{\mathbb{P}^1(u)} C \rightarrow C$. Fix $(\mathrm{id}, \mathrm{id}, \mathrm{id}): C \rightarrow \cc{E}_C \times_{\mathbb{P}^1(u)} C$ as the zero section. The involution $\iota$ on $C$ with respect to the double cover $\varphi_C: C\rightarrow B$ gives us another section for the fibration $ \cc{E}_C \times_{\mathbb{P}^1(u)} C \rightarrow C$, namely, $ (\iota, \mathrm{id}, \mathrm{id})$. It intersects the zero section on the points corresponding to the ramification points $a$ and $b$ of the morphism $C\rightarrow B$. The intersection $(Q_i, Q_i, Q_i)$ where $Q_i= \varphi_C^{-1}(t_i)$, is a singular point on the fibre where it is located if and only if $t_i$ is also a ramification point for $\varphi_{D_{t_i}}: D_{t_i} \rightarrow B$ (see \cite[Cor. 3.2.7, pp.108]{EGA1}). We have two possibilities:

(i) The intersection point is not a ramification point for $\varphi_{D_{t_i}}: D_{t_i} \rightarrow B$.

  Since torsion sections do not meet at a non-singular point (see for
  example \cite[Lemma1.1]{MiPe}) the section given by $(\iota, \mathrm{id}, \mathrm{id})$ has infinite order. The fibration $\cc{E}_C \times_{\mathbb{P}^1(u)} C \rightarrow C$ has positive generic rank and thus, by N\'eron-Silverman's Specialization Theorem, infinitely many fibres with positive rank. Hence the fibration $\cc{E}_C \rightarrow \mathbb{P}^1(u)$ has also infinitely many fibres with positive rank.

(ii) The intersection point is a ramification point for $\varphi_{D_{t_i}}: D_{t_i} \rightarrow B$.

  In this case the curve $C\times_ B D_{t_i}$ is rational. We must consider two cases:

(iia) The curve $D_{t_i}$ induces an independent new section in $\cc{E}\times_B D_{t_i}$ and we are done.

(iib) The generic rank of $\cc{E} \times_B D_{t_i}$ is equal to the generic rank of $\cc{E}$.

 Choose another curve $C \in \cc{L}_1$ in the beginning of the proof. Since by Corollary \ref{principal} only finitely many $D_u$ do not contribute with an extra section after base change, if (iib) holds for almost all curves $C_t \in \cc{L}_1$ then the set $R = \cup_{t \in \mathbb{P}^1} \{b \in B(k); \phi_t: C_t \rightarrow B \text{ is ramified above }\, b\}$ is finite. By Lemma \ref{ramif} below we conclude that all morphisms $\phi_t: C_t\rightarrow B$ ramify above the same points. In this case we start over by fixing a curve $D \in \cc{L}_2$ such that the surface $\cc{E}_D$ has generic rank strictly larger than the generic rank of $\cc{E}$. The curves $C_t \times_B D$ where $C_t$ varies in $\cc{L}_1$ induce an elliptic fibration on $\cc{E}_D$. Since infinitely many of them contribute with a new independent section we will be either in case (i) or case (iia). 

\begin{lem}\label{ramif}
Let $\pi: X \rightarrow B$ be a fibration on curves from a smooth proper surface $X$ defined over a number field $k$ to a smooth proper curve $B$. Let $f: X \rightarrow \mathbb{P}^1$ be a genus zero fibration on $X$ such that the fibres of $f$ are not fibres of $\pi$. Let $\pi_t$ be the restriction of $\pi$ to the fibre $f^{-1}(t)$, for $t \in \mathbb{P}^1(k)$. Let $R= \cup_{t \in \mathbb{P}^1} \{b \in B(k); \pi_t \text{ is ramified above }\, b\}$. Then $R$ is either infinite or equal to $\{b \in B(k); \pi_{t_0} \text{ is ramified above } \,b \}$ for any $t_0 \in B(k)$.
\end{lem}

\subsection{Proof of Theorems \ref{th1} and \ref{th3}}\label{prf}
\textbf{} 
\\
Let $\cc{E}$ be as in the hypothesis of Theorem \ref{th1} or Theorem \ref{th3}. Let $\cc{L}'_1= \{C_t\}_{\{t \in \mathbb{P}^1 \}}$ and $\cc{L}'_2= \{D_u\}_{\{u \in \mathbb{P}^1\}}$ be the two pencils of rational curves constructed in the previous subsections according to the possible minimal models of $\cc{E}$. 
By Corollary \ref{principal} all but finitely many curves $C_t \in \cc{L}_1$ induce a new section in $\cc{E}_{C_t}$ independent of the old sections. For each $t \in \mathbb{P}^1_k$ the pencil $\cc{L}'_{2,t}= \{D_u \times_B C_t\}$ of curves in $\cc{E}_{C_t}$ also satisfies Corollary \ref{principal} and thus for all but finitely many $u\in \mathbb{P}^1_k$ the curve $D_u \times_B C_t$ induces a new section in $\cc{E}_{D_u\times_B C_t}$ independent of the old sections coming from $\cc{E}_{C_t}$.  Thus, after excluding finitely many $t\in \mathbb{P}^1$ and finitely many $u\in\mathbb{P}^1$,
 the surface $\cc{E}_{D_u\times_B C_t}$ satisfies: \[\mathrm{rk}(\cc{E}_{D_u\times_B C_t}(k(D_u\times_B C_t)))\geq \mathrm{rk}(\cc{E}_{C_t}(k(C_t)))+1 \geq \mathrm{rk}(\cc{E}(k(B)))+2.\]  Moreover, since we excluded only finitely many of each $t$ and $u$, by Theorem \ref{posi} we may choose the curve $D_u\times_B C_t$ such that it has infinitely many $k$-rational points.

  \begin{rem}
\emph{The proof of Theorem \ref{th3} breaks down if the $k$-minimal model considered, $X$, is a surface of degree 6 such that the blow up locus of $f$ contains a Galois orbit with 5 points, or if it has degree 3 or 2. Although we are still able to construct a (family of) rational curve(s) defined over the ground field, the generic member of such a family has two connected components, which cannot be used to finish the proof of the theorem. Moreover, $X$ has no families of irreducible conics (i.e., rational curves intersecting the anti-canonical divisor of the rational elliptic surface with multiplicity two) defined over the ground field.}
\end{rem}

\section{Corollaries from analytic number theory}

 To an elliptic curve $E$ over a number field $k$ we can associate a
 sign $W(E|k)$ intrinsically via the product of local signs $W_v(E|k)$
 (for a complete definition of the local sign see for example \cite{Ro}).
 
  The parity conjecture may be stated in the following form:
  
  \begin{conj}\label{parity} Let $E$ be an elliptic curve over a number field $k$. Let $r$ be the rank of its Mordell Weil group. Then \[ W(E|k)=(-1)^r.\]
  \end{conj}
  
\begin{rem}
\emph{The previous conjecture is a weak version of the Birch Swinnerton-Dyer conjecture.}
\end{rem}
  
  Over the field of rational numbers $\QQ$ we know since Wiles that
  $E$ is modular and that $W(E|k)$ is the sign of the functional
  equation of the $L$-function $L(E,s)$. 

   Let $\pi: \cc{E} \rightarrow B$ be an elliptic surface over $k$
   and $U \subseteq B$ an affine open subset over which $\cc{E}$ is an abelian
   scheme. Note \[ U_{\pm}(k)= \{ t \in U(k); W(\cc{E}_t)= \pm 1 \}.\]
   Modulo the parity conjecture we have as well \[ U_{+}(k)= \{ t \in U(k); \text{ rank } \cc{E}_t(k) \text{ is even } \} \text{ and } U_{-}(k)=\{t \in U(k); \text{ rank }\cc{E}_t(k) \text{ is odd }\}.\]

   There are examples for which $W(E_t)$ is constant, but they all
   correspond to isotrivial surfaces (see for example \cite{CS}). In
   the non-isotrivial case and for $B \simeq \mathbb{P}^1$, H. Helfgott
   \cite{He} has shown, under classical
   conjectures, that the sets $U_{\pm}(k)$ are infinite. This is
   established unconditionally in some interesting cases by Helfgott \cite{He},
   \cite{He2} and Manduchi \cite{Ma}. Much less has been done in the
   case $B$ is a genus one curve such that $B(k)$ is infinite. But,
   from previous work cited above, it
   seems reasonable to conjecture the following:

\begin{conj}\label{conj1}
Let $\cc{E} \rightarrow B$ be a non-isotrivial elliptic surface
defined over a number field $k$ such
that $g(B)= 1$ and $B(k)$ is infinite. Then $U_{+}(k)$ and $U_{-}(k)$
are infinite. 
\end{conj}

This allows us to state the following better but conditional result.
\begin{thm}[modulo Conjectures \ref{conj1} and \ref{parity}] 
Let $\cc{E}\rightarrow B$ be a rational elliptic surface satisfying
the hypothesis of Theorem \ref{th1}. Then \[\# \{t \in B(k); r_t\geq
r+3 \}= \infty,\]
where $r$ is the generic rank and $r_t$ the rank of the fibre $\pi^{-1}(t)$.
\end{thm}

\section*{Appendix}
In this appendix we deal with rational elliptic surfaces that have one non-reduced fibre.

 If the blown up points are in non-general position, then the constructions given in Subsection \ref{constructions} may yield reducible curves. Nevertheless, we are still able to deal with some of these cases since some of the special Galois invariant configurations of the base points of a pencil of cubic curves in the plane yield elliptic surfaces with fibre types that are easier to treat, namely, non-reduced fibres. 

 If the surface has a unique non-reduced fibre, then, depending on the structure of a pencil inducing $\cc{E}$, we will be able to prove the rank jumps for infinitely many fibres by first base-changing by a curve in $\cc{L}'_1$ where $\cc{L}_1$ is the pencil of lines through $p_1$ constructed in i) of Subsection \ref{constructions}. The proposition below tells us that the resulting base-changed elliptic surface is still rational and satisfies the hypothesis of Theorem \ref{th2}, i.e., it is a $k$-unirational elliptic surface.

\begin{pro}\label{ng2}
Let $\cc{E}\rightarrow B$ be a rational elliptic surface defined over a number field $k$ such that there is a $k$-birational morphism $\cc{E} \rightarrow \mathbb{P}^2$ contracting the zero section to a point $p_1\in \mathbb{P}^2(k)$. Suppose  $\cc{E}$ has a unique fibre of type $*$, induced by a cubic curve of the form $3m$ or $m \cup 2l$ where $m$ is a line through $p_1$ and $l$ is another line. 

Then for all but finitely many $L_1 \in \cc{L}_1$ the morphism $L'_1 \rightarrow B$ is ramified over the place corresponding to the fibre of type $*$. Moreover, the surface $\cc{E}\times_{B}L'_1$ is $k$-unirational for all but finitely many $L_1' \in \cc{L}'_1$. 
\end{pro}

\noindent\textbf{Proof:} Let $F$ be the non-reduced fibre of $\cc{E}$ given in the hypothesis. 

Suppose first that $F$ is induced by the triple line $3m$ where $m$ is a line through $p_1$. Note that $p_1$ is a base point with multiplicity at least $3$ and thus $f$ factors through the blow up of $p_1$ and two infinitely near points to it. The first blow up of $p_1$ transforms $3m$ into $3m'+2E_1$ where $E_1$ is the exceptional curve above $p_1$ and $m'$ is the strict transform of $m$. The strict transform of $L_1$ intersects $2E_1$, but does not intersect the curve $3m'$. The second blow up is that of $p'_1$, the intersection point of $3m'$ and $2E_1$. Since the strict transform of $L_1$ by the first blow up does not pass through $p'_1$, this curve or its intersection with other divisors is unaffected by the remaining blow ups. Thus the strict transform of $L_1$ by $f$, i.e., $L'_1$ intersects the multiplicity two component of $F$ in a single point. This assures that the map $\varphi_{L_1}:L'_1 \rightarrow B$ is ramified above the place corresponding to $F$.

Now suppose that $F$ is induced by $m \cup 2l$ where l is another line. We may suppose that $L_1$ does not pass through other base points and hence $L_1$ intersects $2l$ at a point that is not in the blow up locus of $f$. This implies that $L'_1$ intersects $F$ at the double component corresponding to the strict transform of $2l$ and thus $\varphi_{L_1}$ ramifies above the place corresponding to $F$.

By Lemma \ref{stillrational}, $\cc{E} \times_B L'_1$ is rational. The existence of a non-reduced fibre in $\cc{E}$ assures that $\mathrm{rank}(\mathrm{Pic}(\cc{E})|k) \geq 4$, since the components of the fibre that do not intersect the zero section contribute with at least two divisors to the Picard group over $k$. By Corollary \ref{principal}, the rank, over $k$, of the surface $\cc{E}\times_B L'_1$ is strictly larger than that of $\cc{E}$, for all but finitely many $L_1 \in \cc{L}_1$, and thus $\mathrm{rank}(\mathrm{Pic}(\cc{E}\times_B L'_1)|k) \geq 5$. By Lemma \ref{unirational}, $\cc{E}\times_B L'_1$ is $k$-unirational. \qed

We apply Theorem \ref{th2} to the surfaces satisfying the hypothesis of the previous proposition. This gives us the following theorem.
 
\begin{thm} Let $\cc{E} \rightarrow B$ be a rational elliptic surface as in Proposition \ref{ng2}. Then there is a finite covering $C \rightarrow B$ such that $C \simeq_k \mathbb{P}^1$ and the surface $\cc{E}_C= \cc{E}\times_B C$ satisfies: 
\[ \mathrm{rank}(\cc{E}_C(k(C)))\geq \mathrm{rank}(\cc{E}(k(B)))+2.\] 

\end{thm}

 As before, we get the following corollary.
\begin{cor}
Let $\pi: \cc{E} \rightarrow B$ be an elliptic surface as in Proposition \ref{ng2}. For $t\in B(k)$, let $r_t$ be the rank of the fibre above the point $t$ and $r$ the generic rank. Then \[ \# \{ t \in B(k); r_t \geq r+2\}= \infty .\]
\end{cor}

\subsection*{Acknowledgements}
\label{sec:acknowledgements}
The author would like to thank Marc Hindry for suggesting this beautiful problem and for several fruitful discussions on the topic, Damiano Testa, Daniel Bertrand, Douglas Ulmer, Jean-Louis Colliot-Th\'el\`ene and Ronald van Luijk for discussions and suggestions that improved the quality of the paper. Finally the author would like to thank the referee for the careful reading of the article as well as for important suggestions that improved the readability and the quality of the paper.

\bibliographystyle{alpha}

\end{document}